\documentclass{amsart}
\usepackage{graphicx}
\usepackage{hyperref}
\usepackage{xcolor}

\newtheorem{theorem}{Theorem}
\newtheorem{lemma}[theorem]{Lemma}

\theoremstyle{definition}
\newtheorem{ques}{Question}

\theoremstyle{remark}

\newcommand{\Z}{\mathbb Z}
\newcommand{\R}{\mathbb R}

\newcommand{\bibtitle}[1]{\emph{#1}}
\newcommand{\dfn}[1]{\textbf{#1}}
\newcommand{\ty}{\nabla\mathrm{Y}}
\newcommand{\yt}{\mathrm{Y}\nabla}

\newcommand{\lk}{\mathrm{lk}}

\newcommand{\modtwo} {\mathrm{mod} \; 2}

\newcommand{\idfour}{I$D_4$}
\newcommand{\iatwo}{I$a_2$}

\title{A brief survey on intrinsically knotted and linked graphs}
\author{Ramin Naimi}
\address{Department of Mathematics,
Occidental College,
Los Angeles, CA 90041, USA.}

\date \today

\begin{document}

%\begin{abstract}    % type your abstract below
%This is a brief, expository survey of some results and open problems
%on intrinsically knotted and intrinsically linked graphs.
%
%\end{abstract}

\maketitle

%%%%%%%%%%%%%%%%%%%%   Start of main body of article

\section{Introduction}

\begin{figure}[ht]

 \centering
 \includegraphics[width=70mm]{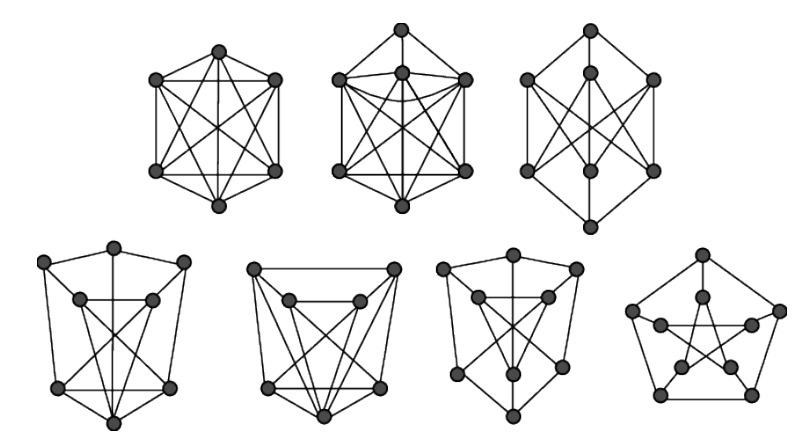}
 \hspace{6mm}
  \includegraphics[width=35mm]{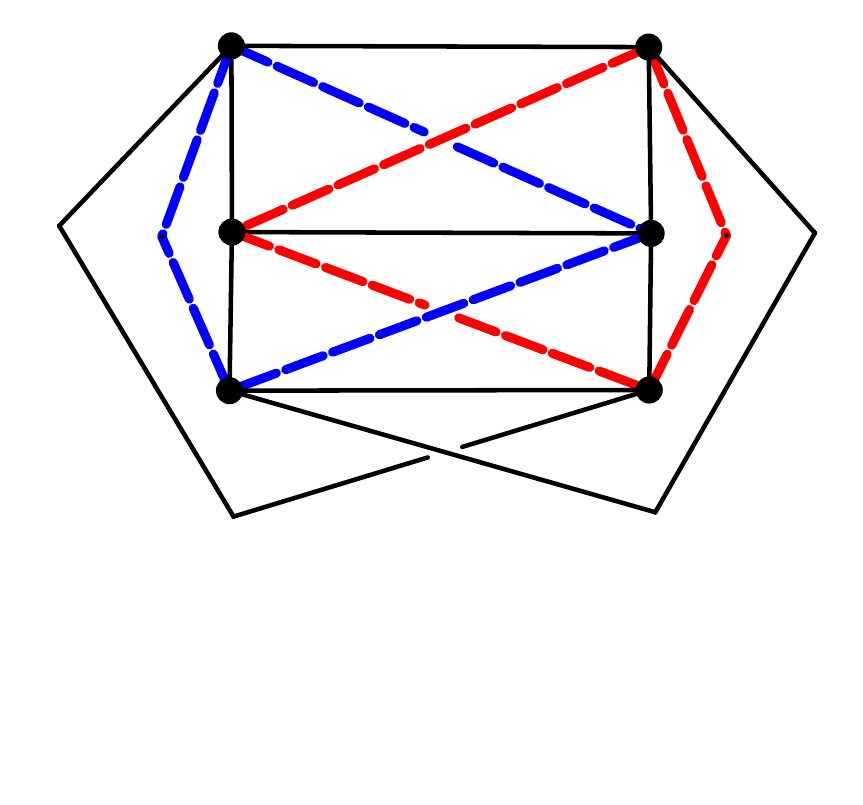}
 \caption{Left: The Petersen family graphs~\cite{WikiPeter}. Right: An embedding of $K_6$, with a nontrivial link highlighted.
 \label {fig-PetersenFamily}
 }
\end{figure}

In the early 1980's, Sachs \cite{sa, sa2}
showed that if $G$ is one of the seven graphs in Figure~\ref{fig-PetersenFamily},
known as the Petersen family graphs,
then every spatial embedding of $G$, i.e.\ embedding of $G$ in $S^3$ or $\R^3$,
contains a nontrivial link --- specifically,
two cycles that have odd linking number.
Henceforth, spatial embedding will be shortened to embedding;
and we will  not distinguish between an embedding and its image.
A graph is  \dfn{intrinsically linked} (IL) if
every  embedding of it contains a nontrivial link.
For example,  Figure~\ref{fig-PetersenFamily} shows a specific embedding of 
the first graph in the Petersen family,
$K_6$, the complete graph on six vertices,
with a nontrivial 2-component link highlighted.
At about the same time, 
Conway and Gordon~\cite{cg}
also showed that $K_6$ is IL.
They further showed that $K_7$ in
\dfn{intrinsically knotted} (IK), 
i.e.\ every spatial embedding of it contains a nontrivial knot.

A graph $H$ is a \dfn{minor} of another graph $G$
if $H$ can be obtained from a subgraph of $G$
by contracting zero or more edges.
It's not difficult to see that if $G$ has a linkless (resp.\ knotless) embedding,
i.e., $G$ is not IL (IK), then every minor of $G$  has a linkless (knotless) embedding \cite{FelLan, NesTho}.
So we say the property of having a linkless (knotless) embedding is \dfn{minor closed} (also called \emph{hereditary}).

A graph $G$ with a given property
is said to be \dfn{minor minimal} with respect to that property
if no minor of $G$, other than $G$ itself, has that property.
Kuratowski \cite{Kur} and Wagner \cite{Wag} showed that
a graph is nonplanar, i.e.\ cannot be embedded in $\R^2$,
iff it contains
$K_5$ or $K_{3,3}$ (the complete bipartite graph on $3+3$ vertices) as a minor.
Equivalently, these two graphs are the only minor minimal nonplanar graphs.
It's easy to check that each of the seven graphs in the Petersen family is minor minimal IL (MMIL),
i.e., if any edge is deleted or contracted, the resulting graph has a linkless embedding.
Sachs conjectured that these are the only MMIL graphs,
which was later proved by
Robertson, Seymour, and Thomas:

\begin{theorem} 
\cite{RST}
\label{thm-rst}
A graph is IL iff it contains a Petersen family graph as a minor.
\end{theorem}

This gives us an algorithm for deciding whether any given graph $G$ is IL:
 check whether  one of the Petersen family graphs is a minor of $G$.

In contrast, 
finding all minor minimal IK (MMIK) graphs has turned out  more difficult.
Robertson and Seymour's Graph Minor Theorem~\cite{RS}
says that in any infinite set of (finite) graphs,
at least one is a minor of another.
It follows that for any property whatsoever (minor closed or not),
there are only finitely many graphs that are minor minimal with respect to that property.
In particular, there are only finitely many MMIK graphs.
If we knew the finite set of all MMIK graphs,
we would be able to decide whether or not any given graph is IK.
So far there are at least 264 known MMIK graphs \cite{FMMNN},
and, for all we know, this could be just the tip of the iceberg ---
we don't even have an upper bound on the number of MMIK graphs.

For $n \ge 3$, a graph is \dfn{intrinsically $n$-linked} (I$n$L)
if every spatial embedding of it 
contains a nonsplit $n$-link (a link with $n$ components).
It was shown in \cite{FNP}
that $K_{10}$ is I3L;
and it was shown in \cite{BowFoi} that 
removing from $K_{10}$ four edges that share one vertex, 
or two nonadjacent edges, yields  I3L graphs; 
but it's not known if they are MMI3L.
Examples of  MMI$n$L graphs
were given
for every $n \ge 3$ in \cite{FFNP}.

Other ``measures of complexity'' have also been studied.
For example:
given any pair of positive integers $\lambda$ and $n\ge 2$,
every embedding of a sufficiently large complete graph 
contains a 2-link with linking number at least $\lambda$ in magnitude \cite{Fla, ShirTan};
contains a nonsplit $n$-link with all linking numbers even \cite{FleDie};
contains a knot $K$ such that 
the magnitude of the second coefficient of its Conway polynomial, i.e.\ $|a_2(K)|$,
%(defined in Section~\ref{sec-IK}),
is larger than $\lambda$ \cite{Fla, ShirTan};
and
contains a nonsplit  $n$-link $L$ such that
for every component $C$ of $L$,  $|a_2(C)| > \lambda$
and for any two components $C$ and $C'$ of $L$,  $|\lk(C,C')| > \lambda$ \cite{FMN}.
No minor minimal graphs with respect to any of these properties are known.

In the following sections we discuss the above, and a few other topics,
in greater detail.

\section{IL graphs}
The proof that $K_6$ is IL is short and beautiful:
There are $20$ triangles (3-cycles) in $K_6$.
For each triangle, there is exactly one triangle disjoint from it.
Thus there are exactly 10 pairs of disjoint triangles, i.e.\ 2-links, in $K_6$.
Any pair of disjoint edges is contained in exactly two such links.
So, given any embedding of $K_6$, 
any crossing change between any two disjoint edges
affects the linking number of exactly two links, and 
the magnitude of each of their linking numbers changes by 1.
Thus, the sum of all linking numbers does not change parity under any crossing change.
Now,  in the embedding of $K_6$ shown in Figure~\ref{fig-PetersenFamily}, 
the sum of all ten linking numbers is odd.
Therefore the same is true in every embedding of $K_6$,
since any embedding can be obtained from any other embedding by isotopy and  crossing changes.
Hence every embedding contains at least one link with odd linking number.

\begin{figure}[ht]

 \centering
 \includegraphics[width=70mm]{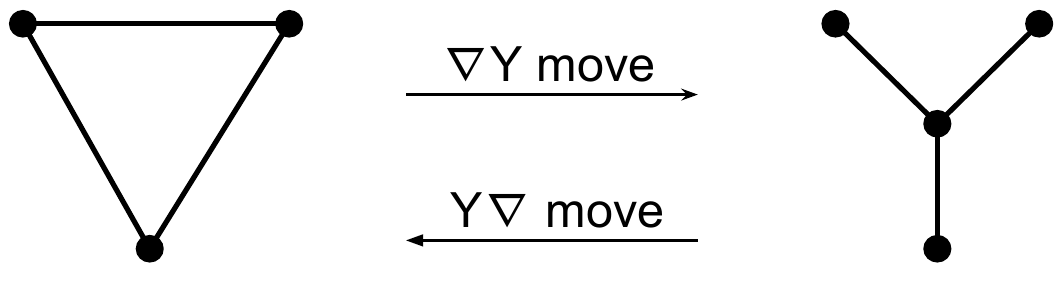}
 \caption{$\ty$ and $\yt$ moves.
% (Other commonly used notation: 
%$\nabla$ =  $\Delta$ = Delta =  triangle; 
%$\mathrm{Y}$ = wye = star;
%move = exchange = transformation = operation.)
 \label {fig-TYmove}
 }
\end{figure}

Sachs observed that a similar argument can be used
to show all seven graphs in the Petersen family are IL.
He also observed that each of these graphs can be obtained from any other
by one or more $\ty$ and $\yt$ moves, as defined in Figure~\ref{fig-TYmove}.
Furthermore, this family is closed under $\ty$ and $\yt$ moves.
It follows, by Theorem~\ref{thm-rst}, that  $\ty$ and $\yt$ moves  preserve the property of being MMIL.

In fact, Sachs observed that a $\ty$ move on any IL graph yields an IL graph;
or, equivalently, that a $\yt$ move on a ``linklessly embeddable'' graph 
yields a linklessly embeddable graph.
The proof of this is elementary and straightforward,
which we outline here.
Suppose $G'$ is obtained from a linklessly embeddable graph $G$ by a $\yt$ move.
Take a linkless embedding $\Gamma$ of $G$,
and replace the Y involved in the $\yt$ move
by a triangle $\nabla$ whose edges are ``close and parallel'' to the edges of the Y.
This gives an embedding $\Gamma'$ of $G'$.
It's easy to show that 
any link in $\Gamma'$ that doesn't have $\nabla$ as a component is isotopic
to a link in $\Gamma$, and hence is trivial.
And any link in $\Gamma'$ that does have $\nabla$ as a component 
is also trivial since $\nabla$ bounds a disk 
whose interior is disjoint from $\Gamma'$.

%which we'll denote as $G \stackrel{\ty}{\to} G'$.

It is also true that a $\yt$ move on any IL graph yields an IL graph,
but the only known proofs of it rely on Theorem~\ref{thm-rst} or the following result of \cite{RST}:
If $G$ has a linkless embedding, then it has a paneled embedding,
i.e., an embedding $\Gamma$ such that 
every cycle in $\Gamma$ bounds a disk
whose interior is disjoint from $\Gamma$.
%together with the following observation
%(whose proof is simple and left to the reader):
%If $G'$ is obtained by a $\ty$ or $\yt$ move on a graph $G$
%and if $H$ is a minor of $G$ such that 
%at least one edge of the $\nabla$ or the $\mathrm{Y}$ that is being replaced in the move
%is deleted or contracted in $G$ to obtain $G'$, 
%then $H$ is a minor of $G'$.

Let's say a graph is $\Z_2$-IL
if every embedding of it contains a 2-link
with linking number nonzero $\modtwo$.
Thus, each of the Petersen family graphs is $\Z_2$-IL.
This, together with Theorem~\ref{thm-rst},
implies that $G$ is IL iff it is $\Z_2$-IL.

It is possible to determine if a graph $G$ is $\Z_2$-IL
by simply solving a system of linear equations, without even using Theorem~\ref{thm-rst}.
We give an outline here.
First, pick an arbitrary embedding $\Gamma$ of $G$,
and compute the linking numbers $\modtwo$ 
for all 2-links (pairs of disjoint cycles) in $\Gamma$.
An arbitrary embedding $\Gamma'$ of $G$ can be obtained from $\Gamma$
by adding some number of full twists between each pair of disjoint edges, plus isotopy
(adding twists is equivalent to letting edges ``pass through'' each other).
Say there are $d$ pairs of disjoint edges in $G$.
Let $x_1, \cdots, x_d$ be variables representing the number of full twists 
to be added to the $d$ disjoint pairs of edges
to obtain $\Gamma'$ from $\Gamma$.
Then the linking number of any 2-link in $\Gamma'$ can be written in terms of
$x_1, \cdots, x_d$  and the linking number of that 2-link in $\Gamma$.
Setting each of these expressions equal to zero
gives us a system of linear equations in $d$ variables.
This system of equations has a solution in $\Z_2$ iff
$G$ is not $\Z_2$-IL.
%An implementation of this algorithm as a Mathematica program is available freely at \cite{??}.
%But for graphs with more than about 30 edges the program become is slow to be of much practical use.
Note that the number of cycles in a graph can grow exponentially with the graph's size,
so this algorithm is exponential in time and space.
In \cite{KKM, vdH}, polynomial time algorithms are given for finding linkless embeddings of graphs.

\section{IK graphs}
\label{sec-IK}
Essentially the same argument that shows the $\ty$ move preserves ILness
also shows the $\ty$ move preserves IKness.
However, the  $\yt$ move does not necessarily preserve IKness \cite{FN}.
For example, there are twenty graphs that can be obtained from $K_7$
by zero or more $\ty$ and $\yt$ moves.
Six of these graphs cannot be obtained from $K_7$ by $\ty$ moves only
--- they require $\yt$ moves also.
And it turns out all these six graphs have knotless embeddings \cite{FN, GMN,HNTY}.

Given two disjoint graphs $G_1$ and $G_2$, 
let $G_1 * G_2$, the \dfn{cone} of $G_1$ with $G_2$,
be the graph obtained  
by adding all edges
from vertices of $G_1$ to vertices of $G_2$,
i.e.,  $G_1 * G_2 = G_1 \cup G_2 \cup \{v_1v_2 \; | \; v_i \in V(G_i)\}$.

For about twenty years, 
the only known IK graphs were $K_7$ 
and its \dfn{descendants}, 
i.e., graphs obtained from $K_7$ by $\ty$ moves only.
It was suspected that $K_{3,3,1,1}$
(the complete 4-partite graph on $3 + 3 + 1 + 1$ vertices)
is also IK.
Recall that $K_5$ and $K_{3,3}$ are minor minimal nonplanar.
Coning with one vertex on each of these graphs 
gives $K_6$ and $K_{3,3,1}$, both of which are in the Petersen family and hence MMIL.
Coning again gives $K_7$ and $K_{3,3,1,1}$;
and $K_7$ was shown to be MMIK;
so it was natural to ask if $K_{3,3,1,1}$ is too.
Foisy~\cite{Foi1} proved that $K_{3,3,1,1}$ indeed is IK.
His technique, partially outlined below, led to finding many more MMIK graphs
later on \cite{Foi2, Foi3, GMN}.

\begin{figure}[ht]

 \centering
 \includegraphics[width=35mm]{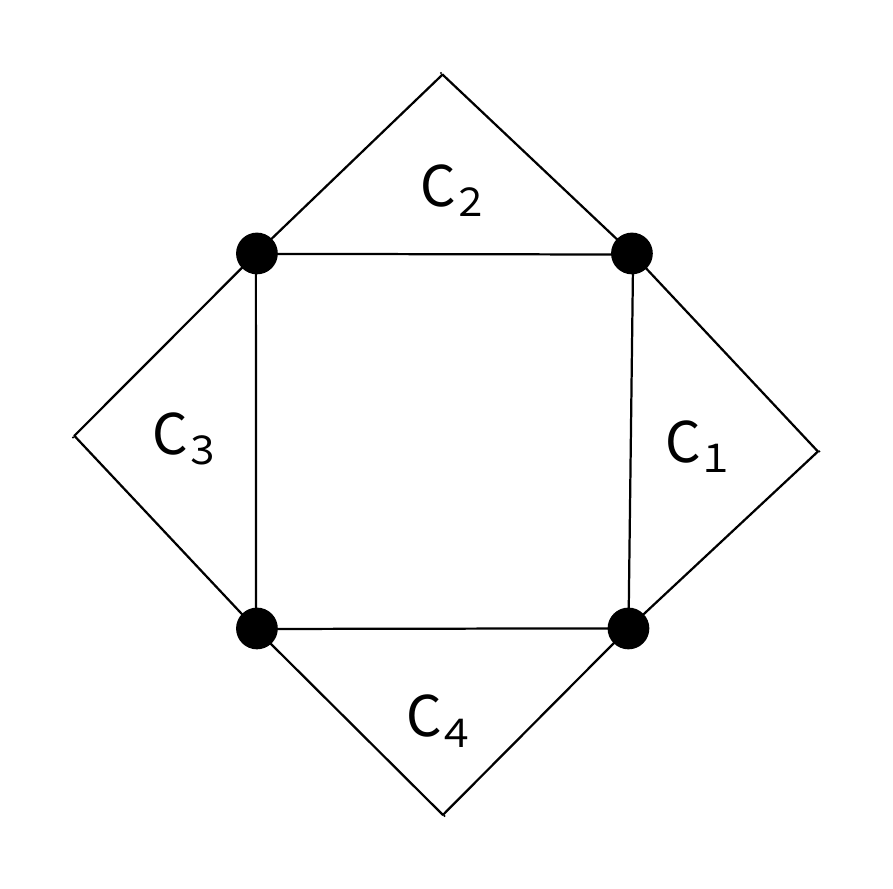}
 \caption{The $D_4$ graph.
 \label {fig-D4-labeled}
 }
\end{figure}

Figure~\ref{fig-D4-labeled} shows a multi-graph (i.e.\ double edges and loops are allowed)
commonly called $D_4$,
with four of its cycles labeled $C_1, \cdots, C_4$.
Let's say an embedding of $D_4$ is \dfn{double linked} $\modtwo$
if  $\lk(C_1,C_3)$ and $\lk(C_2,C_4)$ are both nonzero $\modtwo$.
To show $K_{3,3,1,1}$ is IK, Foisy  proved the following key lemma.
A more general version of the lemma was proved, independently, by Taniyama and Yasuhara \cite{TanYas}.

\begin{lemma}[$D_4$ Lemma]
\cite{Foi1,TanYas}
Every embedding of $D_4$ that is double linked  $\modtwo$
contains a knot $K$ with $a_2(K) \ne 0$  $\modtwo$.
\end{lemma}

Foisy proved that $K_{3,3,1,1}$ is IK 
by showing that every embedding of it
contains as a minor a $D_4$ that is double linked $\modtwo$.
%Subsequently, the $D_4$ Lemma was used in \cite{Foi2, Foi3, GMN}
%to prove other graphs were MMIK.
%There are at least 264 known MMIK graphs \cite{FMMNN}.

Let's say a graph $G$  is \idfour \   $\modtwo$
if every embedding of $G$ contains as a minor a $D_4$ that is  double linked $\modtwo$;
and $G$ is {\iatwo}  $\modtwo$
if every embedding of $G$ contains a knot $K$
such that $a_2(K) \ne 0$ $\modtwo$.
Thus, the $D_4$ Lemma says if $G$ is \idfour\ $\modtwo$ then it's \iatwo\ $\modtwo$.
We also know that if $G$ is \iatwo \ $\modtwo$ then it's IK.
Let's abbreviate these two implications as: \idfour\ $\modtwo$ $\implies$ \iatwo\ $\modtwo$ $\implies$ IK.
It is natural to ask if the converse of each of these implications is also true.

\begin{ques}
\label{Q-Ia2ID4}
{(a)}  IK $\implies$ \iatwo \ $\modtwo$?
\ {(b)} \iatwo \ $\implies$ \idfour \ $\modtwo$?
\ {(c)} IK \ $\implies$ \idfour \ $\modtwo$?

\end{ques}

The  question
``I$n$L $\implies$ I$n$L $\modtwo$?''
%ChdAftRef
is also open.
 
It turns out that every known MMIK graph\footnote{That $K_7$ and 
its descendants are \idfour \  is not in the literature
but is believed to be true if the computer program of \cite{mn} is correct.} 
is  \iatwo \ $\modtwo$ and \idfour \ $\modtwo$.
But this is not necessarily evidence that 
the answer to either part of Question~\ref{Q-Ia2ID4} is yes,
because most of the known MMIK graphs
were found by looking for graphs that are \idfour \ $\modtwo$.

Determining if a graph is \idfour\ $\modtwo$
can be done by solving systems of linear equations \cite{mn}.
If it is true that IK $\iff$ \idfour \ $\modtwo$,
then the algorithm of \cite{mn} can be used to decide whether an arbitrary graph is IK.

%The question
% ``IK  $\implies$ \iatwo?''
% can be thought of as the ``knot version''
% of the statement  IL  $\implies$ $\Z_2$-IL.
% This resemblance  is not quite as superficial as it may seem,
%in light of the $D_4$ Lemma
%and the fact that its proof relies on
%the definition of $a_2$  in terms of linking number.

There may be (a lot) more MMIK graphs than have been found so far;
and trying to find some of them might not be too hard.
For example, one can start with a non-MMIK graph $G$ 
obtained by a $\yt$ move from a known MMIK graph,
and keep adding new edges to $G$
or expanding  vertices of $G$ into edges (the reverse of contracting edges)
until one obtains an IK graph.
But finding more and more MMIK graphs
doesn't seem to have advanced our understanding of IK graphs very much.
In trying to understand IK graphs better,
another approach has been to try to classify
all IK graphs with a given number of edges.
For example, it has been shown that there are no IK graphs with 20 or fewer edges,
and the only IK graphs with exactly 21 edges are $K_7$ and its descendants \cite{BarMat, JKM, LKLO, Mat};
IK graphs with 22 edges have also been partially classified \cite{KMO, LLMO}.
However, this approach doesn't seem to have led
to significant insights or advances in the theory either.

\section{Miscellaneous facts and open problems}

In \cite{NPS} it was shown that 
if $G$ is IK and $e$ is an edge of a 3-cycle in $G$,
then $G \setminus e$ is IL.
The following related questions might be useful
in trying to answer Question~\ref{Q-Ia2ID4}.

\begin{ques}
Suppose $G$ is IK. 
(a) Is $G \setminus e$, or $G / e$,  IL for every edge, or for some edge, $e$ of $G$?
(b)  Does $G$ have at least two distinct nonsplit links?
\end{ques}

Sachs observed that a graph $G$ is non-planar iff
the graph $G * v$, the cone of $G$ with one vertex $v$,
is IL.
This can be seen as follows.
If $G$ is nonplanar, then
it contains $K_5$ or $K_{3,3}$ as a minor.
So $G*v$ contains $K_5 * v = K_6 $ or $K_{3,3} *v = K_{3,3,1}$ 
as a minor, and hence $G*v$ is IL.
Conversely, if $G$ is planar,
it is easy to construct a linkless embedding of $G*v$:
embed $G$ in the plane, 
put $v$ above the plane, 
and connect $v$ with straight edges to all vertices of $G$.

A graph $G$ is said to be \dfn{$n$-apex} if there exist  vertices $v_1, \cdots, v_n$ in $G$
such $G - \{v_1, \cdots, v_n\}$, 
i.e., the graph obtained by removing $\{v_1, \cdots, v_n\}$ 
and all edges incident to them,
is planar.
A 1-apex graph is  called \dfn{apex}.
Thus, by above, apex graphs are not IL.
It can similarly be shown that
2-apex graphs are not IK.
In fact, $G$ is planar iff
the graph $G*v*w$ (i.e., $G*K_2$) 
is not IK \cite{BBFFHL, OzaTsu}. 
The reason is similar to the one given above:
If $G$ is nonplanar, 
then $G*v*w$ contains $K_5 * v*w = K_7 $ or $K_{3,3} *v*w = K_{3,3,1,1}$ 
as a minor; 
and since both of these graphs are IK,  $G*v*w$ is IK too.
If $G$ is planar, 
we can construct a knotless embedding of $G*v*w$ as follows:
embed $G$ in the plane, 
 put $v$ above the plane, put $w$ below the plane, 
and connect $v$ and $w$ to all vertices of $G$ and to each other with straight edges.
The list of all minor minimal non-$n$-apex graphs is not known for any $n$, even $n=1$.
A detailed survey of results on apex and $2$-apex graphs can be found in \cite{FMMNN}.

The crossing number $C(K)$ of a knot $K$ is the fewest number of crossings among all regular projections of $K$.
It's easy to see that for every $n$, 
the set $\{K \; : \; C(K) \le n \}$ is finite;
so $A(n)= \max\{|a_2(K)| \; : \; C(K) \le n  \}$ is well-defined and finite.
As mentioned before, 
given a fixed $n$,
every embedding of a sufficiently large complete graph 
contains a knot $K$ with $|a_2(K)| > A(n)$;
hence, every embedding of a sufficiently large complete graph 
contains a knot with  crossing number larger than $n$.
The \dfn{bridge number} of a knot $K$ is the minimum number of local maxima
with respect to height ($z$-coordinate in $\R^3$) 
among all isotopic embeddings of $K$.
Given any integer $n \ge 2$,
there are infinitely many $n$-bridge knots. 
So the above argument for crossing number doesn't work for bridge number.
This leads to the question:
Does there exist, for each $n$,
a graph $G$ such that
every embedding of $G$
contains a knot with bridge number at least $n$?

Suppose $G'$ is obtained by a $\ty$ move from $G$.
It turns out that if $G$ has any of the following properties, then $G'$ has  that property too:
IL; IK; \iatwo; \idfour; I$n$L; nonplanar; non-$n$-apex.
The proofs for all of these are elementary and short,
and most of them are similar to the one we saw for IL.
But, curiously, IL is the only property from the above list
known to be preserved by $\yt$ moves.
%There are known counterexamples for most of the others
%(IK, I3L, nonplanar, ???).
%For example,
%a $\yt$ move on $K_{3,3}$, which is nonplanar,
%yields $K_5$ minus one edge, which is planar.

The \dfn{complement} of a graph $G$
is a graph $G^c$ with the same vertices as $G$
and with exactly those edges not in $G$.
In \cite{BHK} it was shown that
if $G$ has 9 or more vertices, then $G$ or $G^c$ is nonplanar.
This result is sharp:
there exists a graph $G$ on 8 vertices
such that both $G$ and $G^c$ are planar.
We can ask a similar question of IL graphs:
What is the smallest integer $v$ such that
for every graph $G$ with $v$ vertices,
$G$ or $G^c$ is IL?
Here is a partial answer.
In \cite{Mad} it was shown that
%ChdAftRef
for all $n \le 5$,
if $G$ has $v$ vertices, $e$ edges,
and $e > nv - {n+1 \choose 2}$,
then $G$ contains $K_{n+2}$ as a minor.
Now,
 $K_{15}$ has 105 edges,
so if $G$ has 15 vertices,
then $G$ or $G^c$ has at least $\lceil 105/2 \rceil = 53$ edges.
Letting $n =4$, we have
$nv - {n+1 \choose 2 } = 4(15) - {5 \choose 2} = 50$;
since $53 > 50$, 
 $G$ or $G^c$ contains $K_6$ as minor,
and hence is IL. But this is not sharp.
In \cite{PavPav}, it is shown that:
(i)~if $G$ has 13 vertices, then $G$ or $G^c$ is IL,
and (ii)~there is a graph $G$ with 10 vertices
such that neither $G$ nor $G^c$ is IL.
The question for graphs with 11 or 12 vertices remains open.

%ChdAftRef
One can similarly show that
for every graph $G$ with 18 or more vertices,
%ChdAftRef
$G$ or $G^c$
contains a $K_7$ minor and hence is
IK.
%and, given any $n$, $n$ disjoint copies of  \iatwo, \idfour, %I$n$L, 
%and non-$n$-apex.
%For each of these properties, 
It is unknown what the minimum number of vertices is
that would guarantee that $G$ or $G^c$ 
is IK.
%has that property.

Let's say a graph is \dfn{strongly intrinsically linked} (SIL)
if every embedding of it contains a 2-link
with linking number at least 2 in magnitude.
Then $K_{10}$ is SIL, since,
by \cite{FNP}, 
$K_{10}$ contains
a 3-link two of whose linking numbers are nonzero,
and by \cite{Fla}, any embedded complete graph that contains
such a 3-link contains a ``strong'' 2-link.
What about $K_9$?
By \cite{FNP},
$K_9$ does not contain such a 3-link;
but it's not known whether or not $K_9$ is SIL.

A \dfn{digraph} (directed graph) is a graph each of whose edges
is oriented.
A \dfn{consistently oriented} cycle in a digraph
is a cycle $x_0, x_1, \cdots, x_n$, where $x_n = x_0$,
such that each edge $x_i  x_{i+1}$ is oriented
from $x_i$ to $x_{i+1}$.
A digraph is said to be I$n$L (resp.\ IK) if every spatial embedding of it
contains a nonsplit $n$-link 
(resp.\ nontrivial knot)
consisting of consistently oriented cycles.
In \cite{FleFoi}, an IK digraph and an I4L  digraph are constructed.
It is not known whether there exists an I$n$L digraph with $n \ge 5$.
In \cite{FoiHowRic} it was shown that
(unlike all the other graph properties we have discussed)
the property of having a linkless embedding is not minor closed for digraphs.

%\bigskip
%
%(There is various terminology for this in the literature: 
%$\nabla$ =  $\Delta$ = Delta =  triangle; 
%$\mathrm{Y}$ = wye = star;
%move = exchange = transformation = operation.)

%\begin{figure}[ht]
%
% \centering
% \includegraphics[width=80mm]{triangleYmove.pdf}
% \caption{$\ty$ and $\yt$ moves.
% \label{triangleYmoveFigure}
% }
%\end{figure}

%\section*{Acknowldedgements}

%%%%%%%%%%%%%%%%%%%%   End of main body of article
%
%                             References
%
%   BiBTeX users uncomment the following line:
%
%\bibliographystyle{gtart}
%

\end{document}